%Paper: alg-geom/9305003
%From: grassi@msri.org (Antonella Grassi)
%Date: Wed, 5 May 93 09:02:05 PDT

%Ams-TeX 2.1
\input amstex
\documentstyle{amsppt}
%Use the lines below instead of magstep and hsize line
\hoffset0.5in\voffset0.6in
\magnification\magstep1\hsize6truein\vsize8.5truein
\rightheadtext\nofrills{Equidimensional models of elliptic threefolds}
\leftheadtext\nofrills{A. Grassi}
\topmatter
\title
Log contractions and equidimensional models of elliptic threefolds
\endtitle
\author A. Grassi \endauthor
\address  Mathematical Sciences
Research Institute, 1000 Centennial Drive, Berkeley, CA 94720.
 \endaddress
\NoBlackBoxes
\thanks
 Research at MSRI supported in part by NSF grant \#DMS 9022140
\endthanks
\email
grassi\@msri.org
\endemail

\endtopmatter

\def\lco{\llcorner}
\def\rco{\lrcorner}
\def \bs{\bar S }
\def\bda{\bold\Delta}
\def\bla{\bold\Lambda}
\def\bsa{\bold\Sigma}
\def\kla{ K_S + \bold\Lambda}
\def\bx{\bar X}
\def \a{\alpha }
\def \g{\Gamma }
\def\so{{S _0}}
\document

\vskip  0.5in

This work was initially motivated by Miranda's work on elliptic  Weierstrass
threefolds.
 An elliptic variety is a complex (irreducible and reduced) projective variety
together with a morphism $\pi : X \rightarrow S$ whose fiber over a general
point is a smooth elliptic curve.  For example, if the elliptic fibration has a
section, $X$ can be locally
expressed as a hypersurface [De].  This is called the Weierstrass form of the
fibration; all such fibrations are (by definition) equidimensional. If
$dim(X)=2$, the fibration is
of course always equidimensional. This is in general not true in higher
dimension (see [U] for a non trivial example).

Miranda [Mi] describes a smooth equidimensional (flat) model for any elliptic
 Weierstrass threefold; such models occur naturally in the study of moduli
spaces.

Each birational class of equivalent elliptic fibrations over $\bs$ is
associated
to a log variety $(\bs, \bla _{\bs})$ (1.9.2)-(1.10).
 Here (in \S 2) we use  Minimal model theory to link birational maps of log
surfaces
(log contractions) to equidimensional  fibrations of elliptic threefolds.
 In particular we give a necessary and sufficient condition for an elliptic
threefold  $X \to S$ to be birationally equivalent to an equidimensional
elliptic fibration $\bx \to S$, where $\bx$ has terminal singularities and $S$
is the original basis of the fibration (Theorem 2.5).
We call this an {\it equidimensional model}.

As a corollary we show that an elliptic fibration of positive Kodaira dimension
has a \it minimal model \rm with an \it equidimensional \rm birationally
equivalent elliptic fibration
 satisfying certain good properties (Corollary 2.7) that naturally generalize
the case of elliptic fibrations of surfaces.

 Miranda derives his result analyzing the local equation
of the  Weierstrass model.  He shows, checking case by case, that after blowing
up the surface a sufficient number of times it is possible
to resolve the singularities of the threefold preserving the flatness. We
present
a summary of Miranda's results in \S 3.

In \S 4 we apply the results in \S 2 to the case of Weierstrass models and give
a global explanation of Miranda's algorithm (see 4.6).
The argument provided is intrinsic, and the proof does not involve any case by
case checking.

\bigpagebreak

 This paper originates from a section of my 1990 Ph.D. thesis at Duke
University.
The bulk of section 4 and of the arguments  of 2.4 and 2.5 are part of my
thesis.  The remaining part of the research was started while traveling in the
Italian Alps in
the Summer of '92 and continued
at MSRI in the Spring 1993. I would like to thank
D. Morrison, who also supervised the beginning of this project, M. Gross,
for reading various versions of this work,
 and
the people at MSRI for providing a serene work environment.

\vskip  0.2in

{\bf \S   1. Some facts }

\medpagebreak

\noindent (1.1) Let $\pi : X \to S$ be any elliptic fibration.
$$    \bsa _{X/S} =\{  P \in S  \ \text{ such that } \ \pi \text{ is not smooth
over } P \}$$
\noindent is called \it  the ramification locus \rm of $\pi$.
  I. Dolgachev [Do] has showed that if $S$ is smooth, then $    \bsa _{X/S}$ is
a divisor.

We write $    \bsa$ for $    \bsa _{X/S}$ when the fibration is well
understood.

Miranda considers elliptic fibration of threefolds in \it Weierstrass \rm form.
Then the standard equation is:
$$y^2  = x^3 + a(s,t)x +  b(s,t), \tag  1.1.1$$  where the elliptic morphism
sends the point
$(x,y,s,t) $ to $(s,t)$.

 The corresponding equation for the ramification locus turns out to be: $$ 4
(a(s,t))^3 + 27(b(s,t))^2 =  0. $$
 If $a$ and $b$ are chosen generically, this fails to have simple normal
crossings.

\bigpagebreak

\noindent (1.2) We can associate to each point $Q \in S -     \bsa$ the \it
$J$-invariant \rm of the smooth elliptic curve $\pi ^{-1}  (Q)$, where $J$ is
the elliptic modular function defined on the complex upper half plane: $\bold
J: S \dasharrow \Bbb P^1$ [Kd].

 If the elliptic threefold is written is the  Weierstrass form (1.1.1), then
 $$\bold J(s,t) = \frac{4 a^3}{4 (a(s,t))^3 + 27(b(s,t))^2}.$$
In general $\bold J$ is not a morphism:

\smallpagebreak

\noindent{\bf Example 1.2.1}
If the elliptic  variety is defined by
$$y^2 = x ^3 + a x + b,$$
then $\bold J = (4 a^3) \big/ (4 a^3 + 27 b^2).$
If $a= s$  and $b=t$ then the $\bold J$ map is not defined at $ s= t =0$.
Note that $    \bsa = (4 s^3 + 27 t^2)$ has a cusp at $(0,0)$.

\bigpagebreak

\definition{Definition 1.3}A morphism $\pi : X \to S$ is equidimensional if
every fiber
 of $\pi$
has dimension equal to $dim(X) - dim(S)$.
\enddefinition

If $\pi :  X \to S$ is a smooth elliptic surface, then $J: S \to \Bbb P^1$ is a
morphism
and $\pi$ is equidimensional [Kd]. This is not the case in higher dimension [U]
and Kodaira's analysis of the elliptic surfaces cannot be extended directly.

\bigpagebreak

\noindent (1.4) Fix a general (smooth) $P \in D_k$ a reduced component of $
\bsa$ and denote by $X_k$ the singular fiber
over the point $P$. By taking $P$ ``general" we restrict ourselves to an
elliptic
 fibration over a small disc, which is the case of elliptic fibration of
surfaces.
 Kodaira has classified all such $X_k$ [Kd; Th. 6.2, p. 1270].

 Consider a small loop around
$P$ and the induced  monodromy transformation (Kodaira shows that $\bold J$ is
holomorphic in a  neighborhood of the point $P$).

The following table [Kd; Table I, p. 1310--Table II, p. 1346] relates the value
of $\bold J$ to the eigenvalues
of the monodromy matrices and the type of singular fiber $X_k$. The two
eigenvalues
of the monodromy matrix  $A_k$ (when finite) are complex conjugates of each
other: we write the one which lies in the closed upper-half plane as $ e^{ 2
\pi  i a_k }$.

 If the monodromy is finite $ 12a_k = \chi(X_k) $ is the Euler characteristic
of the singular fiber. Otherwise  $\bold J(P) = \infty$, in fact the matrix has
finite period if and only if $\bold J \neq \infty$ [Kd].
Note that $\bold J$ has  pole of order $b=\chi(X_k)$ at $P$ if $X_k$ is of type
$_mI_b$
and of order $b=\chi(X_k) -6$ if $X_k$ is of type $I^*_b$.

Set $a_k(_mI_b) =0$ and $a_k (I^*_b)=1/2$. Kodaira shows that if $\pi : X \to
S$ is a
smooth elliptic \it surface \rm, then
$$12\pi _* (K_{X/S}) \simeq \Cal O _S (\sum 12 a_k D_k ) \otimes J_{\infty}
\tag 1.4.1 $$
where $D_k$ and $a_k = a(D_k)$ are defined as above and $J _{\infty}$ is the
divisor of poles of the elliptic modular function $\bold J$.

\newpage

\indent \phantom{kkkkkkkkkkkkkkkkkkkkkkkkk}{\bf TABLE 1:} Kodaira's Singular
Fibers

\noindent $$
\spreadlines{4.5\jot}
\alignat 8
\bold J & \qquad A _k & \qquad \chi (X_k) &  \qquad \text{e-value}&  \qquad
a_k = a(D_k) &  \qquad \text{type} \qquad \phantom{zumpapa} \\
&\\
\bold J(P) = 0 & \qquad
{ \left(\smallmatrix +1 & +1\\
-1  & 0 \endsmallmatrix \right) }
 & \qquad
 2 &  \qquad  e^{ 2 \pi  i 2 /12 } &  \qquad   \frac{1}{6}  &  \qquad II
\qquad \phantom{zumpapa}\\
\bold J(P) =  0  & \qquad
 { \left(\smallmatrix
0  & -1\\
+1  & 1 \endsmallmatrix \right) }
   & \qquad
 10  &  \qquad e^{ 2 \pi  i 10 /12 }   &  \qquad \frac{5}{6}   &  \qquad II^*
\qquad \phantom{zumpapa}\\
\bold J(P) = 0 & \qquad
{ \left(\smallmatrix
-1  & -1\\
+1  & 0 \endsmallmatrix \right) }
 & \qquad
8 &  \qquad e^{ 2 \pi  i 8 /12 }  &  \qquad \frac{2}{3}  &  \qquad IV^* \qquad
\phantom{zumpapa} \\
\bold J(P) = 0 &  \qquad
{ \left(\smallmatrix
0  & +1\\
-1  & -1
\endsmallmatrix \right) }
 & \qquad
 4 &  \qquad e^{ 2 \pi  i 4 /12 }  &  \qquad \frac{1}{3} &  \qquad IV \qquad
\phantom{zumpapa}\\
\bold J(P) = 1 &  \qquad
{ \left(\smallmatrix
0  & +1\\
-1  & 0
\endsmallmatrix \right) }
 & \qquad
3 &  \qquad e^{ 2 \pi  i 3 /12 }  &  \qquad \frac{1}{4} &  \qquad III \qquad
\phantom{zumpapa} \\
\bold J(P) = 1 &  \qquad { \left(\smallmatrix
0  & -1\\
+1  & 0
\endsmallmatrix \right) }
 & \qquad
9  &  \qquad e^{ 2 \pi  i 9 /12 }  &  \qquad \frac{3}{4} &  \qquad III ^*
\qquad \phantom{zumpapa} \\
\text{ regular } &  \qquad
{ \left(\smallmatrix
+1 &  0\\
0  & +1
\endsmallmatrix \right) }
&  1  &  \qquad +1 &  \qquad 0 &  \qquad _mI_0 \\
\text{ regular } &  \qquad  { \left(\smallmatrix
-1 &  0\\
0  & -1
\endsmallmatrix \right) }
 &  6  &  \qquad -1 &  \qquad \frac{1}{2} &  \qquad I ^* _0  \qquad
\phantom{zumpapa} \\
\text{ pole of order } b  &  \qquad
{ \left(\smallmatrix
+1 & b\\
0  & +1
\endsmallmatrix \right) }
 & b &  \qquad &  \qquad  0 &  \qquad _m I_b  \qquad \phantom{zumpapa} \\
\text{ pole of order } b  &  \qquad { \left(\smallmatrix
-1 & -b\\
0  & -1
\endsmallmatrix \right) }
 & b+6 &  \qquad &  \qquad \frac{1}{2} &  \qquad I ^* _b  \qquad
\phantom{zumpapa}
\endalignat
$$

Kawamata  and Fujita have generalized (1.4.1) to higher dimensions.

\definition{Definition 1.5} A divisor $D$ has {\rm simple normal crossings } if
it is the
union of smooth irreducible components intersecting transversely.
\enddefinition

If $X$ and $S$ are smooth and $    \bsa _{X/S} $ has simple normal crossings
$\bold J$ is a morphism and $\pi _* (K_{X/S})$ is an invertible sheaf [Ka].
Note that in the example 1.2.1  $\bold J$ is not a morphism: in fact $    \bsa$
does not have simple normal crossings (see also \S 5).

Assuming  $    \bsa$ to have simple normal crossings is a way to control how
the components of $    \bsa$ intersect.
In particular this assumption is essential for the proof of Kawamata's results,
in extending the variation of Hodge structure over the ramification locus.

\smallpagebreak

Kawamata [Ka] also shows that
$$12\pi _* (K_{X/S}) \simeq \Cal O _S (\sum 12 a_k D_k ) \otimes J_{\infty}
\tag   1.6 $$
where $D_k$ are irreducible components of $    \bsa$ and $a_k = a_k(P)$ are
defined in Table 1 and $P \in D_k$ is a general point.

\smallpagebreak

\noindent (1.6') \phantom{YUP} If we write $ J_{\infty} = \sum b_j B_j$,
then $J$ has a pole of order $b_j$ on the generic point of $B_j$.

  Thus $12  \pi _* (K_ { X/S})$ is a divisor supported on $    \bsa _{X/S}$.
\smallpagebreak

Under the same assumptions, Fujita [F] has proved the following formula for the
canonical bundle of a smooth elliptic threefold:
$$  \  \phantom{mnb} mK_ X  = \pi ^*\{ m K_S + m \pi_* (K_{X  /S}) + m\sum
 \frac {n_i - 1}{ n_i}  Y_i\} + mE -  mG  \tag   1.7 $$
where the fiber over the general point of
$Y_i$ is a multiple fiber of multiplicity $n_i$, $m$ is a multiple of the $ \{
n_i \}$'s; $mE$  and $mG$ are effective divisors and $codim \pi (G) \geq 2$.
Furthermore
$ m\pi^* (\sum
 \frac {n_i - 1}{ n_i}  Y_i) +m E -  mG $ is an effective divisor.
If the fibration has a section, then $n_i =0, \forall i$.

For convenience, we give the following:

\definition{Definition 1.8}
Let $\pi :   X \rightarrow S$ be an elliptic fibration between smooth
varieties. Assume that $    \bsa$ is a divisor with simple normal crossings.
Let:
$$ \aligned \bda  _{X/S} & =  \sum  a_k D_k + (1/12) J _\infty\\
\bla _{X/S} &=  \bda  _{X/S} + \sum \frac{n_i - 1}{ n_i}  Y_i \text{. }
\endaligned
 $$
\enddefinition

\remark{Remarks }

\noindent  (1.9.1)  $\bda$ and $\bla$  are effective $\Bbb Q$-divisors (see for
example,
[Wi]).

\noindent (1.9.2) Since (1.8) and  (1.6') involve only  the generic point
of an irreducible component of the ramification locus, we can define $\bla
_{\bx / \bs }$ for any elliptic fibration between varieties with isolated
singularities.
 Then $\bla _{\bx/ \bs } = \bla _{\bs } $ depends only on the class of
birationally equivalent elliptic fibration over $\bs $.

\smallpagebreak

\noindent (1.9.3) We say that $\bar \pi : \bx \to \bs$ has {\it  no multiple
fiber over divisors,}  if $n_i =0$, in $\bla _ {\bx / \bs}, \ \forall i$.

\smallpagebreak

\noindent (1.9.4) If $ \psi : \bs  \to   \bs _n$ is a birational morphism,
consider
 the induced elliptic fibration $\bx \to  \bs _n$; then
$\bla _{ \bs } = \psi _*(\bla _{   \bs _n})$.
\endremark

\medpagebreak

 In the notation of 1.8 we have

$$ \phantom{mnb}K_ X  \equiv \pi ^* ( \kla _S ) + E -  G, \tag  1.7'$$
where
$=$ indicates equality of divisors or sheaves, while $\equiv$ denotes numerical
equivalence of $\Bbb Q$-divisors.

\bigpagebreak

Note that the coefficient of the irreducible components of  $ \bla $ in (1.8)
are  non negative rational numbers smaller than $1$:  $(S, \bla)$ is thus
a  {\it log variety} (see foe example [Ko et al]).

Then each birational class of equivalent elliptic fibrations over $\bs$, as in
1.9.2, is associated to a log variety $(\bs, \bla _{\bs})$. $\psi$ as in
(1.9.4)
 is a (log)
morphism between the log varieties $(\bs, \bla _{\bs})$ and $(\bs_n, \bla
_{\bs_n}).$

 If $\phi : S \to (\bs, \bla _{\bs} )$ is the blow up of a log variety, then
there is not a unique choice for  $\bla _S$ (\it the birational transform of
$\bla _{\bs}$,
 \rm) so that $\phi$ is a log morphism between
$(S, \bla _S)$ and $(\bs, \bla _{\bs})$ (see [Ko et al; Ch 2]).
When the log variety is associated to an elliptic fibration, the birational
transform of $\bla _{\bs}$ is determined by the pullback of the elliptic
fibration.
 We will see in section 4 that $\bla _S$ is determined by $\bla _{\bs}$
if $S$ is smooth and the fibration as a section. This is not the case if there
are
 multiple fibers over divisors.

\definition{Definition   1.10} Let $(S,\bla)$ be a log surface with log
terminal singularities (see for example [KMM] or [Ko. et. al]).
A birational morphism $\psi : S \rightarrow \so$, with exceptional divisor $\g$
is called a  ($\kla$)-extremal contraction or a \it log extremal \rm
contraction if
$(\kla) \cdot \g <0$.
\enddefinition

 \vskip 0.15in

 If $X$ is any smooth minimal elliptic surface, then $K_X \equiv \pi ^*(\kla)$
and the fibration is obviously equidimensional.  Although this is not the
always case in
 higher dimension (see 1.11 below and more generally 1.14), it is known that
there exists \it some \rm
suitable birational model of the fibration where at least one of the two
properties above is satisfied [G2] and [N2].

 In \S 2 we link log contractions  $(S, \bla _S) \to (\bs, \bla _{\bs})$ to the
existence
of  equidimensional models over $\bs$.

In the last part of this section we  use Mori's algorithm [Mo] to characterize
those elliptic fibrations $X \to S$ which have
 birationally
equivalent fibration $\bx \to S$ satisfying the pullback property.

$\bx$ has at most terminal singularities,
and the canonical bundle is numerically equivalent to the pullback of a
divisor on the surface so that the canonical divisor is $\bar \pi$-nef. ($K
_{\bx}$ is said to be $\bar \pi$-nef when $K _{\bx }\cdot \g \geq 0$
for all curves $\g$ contracted by $\bar \pi$, see for example, [W].)

\medpagebreak

\demo{Example 1.11  (M.Gross)} Take $E \times \Bbb F ^1$, where $E$ is a smooth
elliptic curve and $f: \Bbb F ^1 \to \Bbb P ^1$ the minimal ruled surface
with exceptional  $-1$ curve $\g$. Performing 2 ``suitable" logarithmic
transformation on 2
different fibers of $f$ (say $f_1$ and $f_2$) we obtain a smooth elliptic
threefold $\pi : X \to \Bbb F ^1$
with $\bsa $ the (disjoint) curves $f_1$ and $f_2$.
The type of the general fiber over each $f_i$ is $_3 I_0$. Furthermore
$K_X \equiv \pi ^* (\kla)$.

Consider the elliptic fibration $\epsilon : X \to \Bbb P^2$, induced by  $ \psi
: \Bbb F ^1 \to \Bbb P ^2$  ($\g$ is contracted by $\psi$). Then
$\bla _{\Bbb P ^2} = \frac{2}{3} h +   \frac{2}{3} h $ ($h$ is the class of a
line),
$\bla _ {\Bbb F ^1} = \frac{2}{3} f_1 +   \frac{2}{3} f_2 $ while
$ \psi ^* (\bla _{\Bbb P ^2}) =\frac{2}{3} f_1 +   \frac{2}{3} f_2 +
\frac{4}{3} \g$.

    $X$ cannot have a birational equivalent model $\bar \pi : \bx \to \Bbb P
^2$ such that
$K_{\bx} \equiv \bar \pi ^*(K_{\Bbb P^2} + \bla)$. In fact
 $K_X \equiv \epsilon ^*(\kla) - \frac{1}{3} \g $  (1.12). See also  4.3.2.
\enddemo

\smallpagebreak

 More generally the following hold:

\proclaim{Lemma 1.12} Let $ \pi _i : X_i \to S, \ i= 1,2$ be birationally
equivalent
elliptic fibrations, $X_i$ with terminal singularities, $S$ with log terminal
singularities.
 Assume that $$K_{X_i} \equiv \pi _i ^* (\kla) + D_i, $$ for some $\Bbb
Q$-divisor
$D_i$. Then $D_1$ is effective if and only if $D_2$ is effective.
\endproclaim

\demo{Proof} Let $f_i: X \to X_i$ be a  common resolution with exceptional
divisors $F ^i _k$.
 Then
$$ K_X \equiv f_i ^*(K_{X_i}) + \sum e_k F ^i _k  \equiv \pi ^*(\kla) +
f_i ^*(D_i) +\sum e_k F ^i_k, \ \quad i=1,2$$
 where $\sum e_k F_k  $
is an effective divisor. Then $f_1 ^*(D_1) +\sum e_k F ^1_k  \equiv f_2 ^*(D_2)
+\sum e_k F^2_k$ and $ |m(f_i ^*(D_i) +\sum e_k F ^i_k )| =| m( f_i ^*(D_i)) |
\ \forall m>>0$.
 The statement follows.
 \qed
\enddemo

\smallpagebreak

\proclaim{Theorem   1.13} Let $\pi : X \to S$ be an elliptic fibration, $X$
with terminal
and  $(S, \bla)$ with log terminal singularities.
Assume that $K_X \equiv  \pi ^*  (K_ S +  \bla ) + E-G $, as in (1.7),
with $E-G$ effective. Then

(1.13.1) $ \quad K _  {X} \equiv  \pi^* ( K_{S} + \bla )+ E  . $

(1.13.2) Furthermore there exists a birationally equivalent fibration
 $\bar \pi  : \bx  \rightarrow S $ such that $\bx$ has only $\Bbb Q$-factorial
terminal singularities, $K_{\bx }$ is $\bar \pi $-nef and $$rK_{\bx} = \bar \pi
^* r (K_ S +  \bla ) .$$
\endproclaim

\demo{Proof}
(1.13.1) Let $E -G =D$ be an effective divisor (1.7). Write  $ D=L + M _+ - M
_- $,
with $L, \  M _+, \  M _-$ effective divisors, such that $L$ is the part with
no non-flat component, while
$ M_ +$ and $ M_ -$ map via $\pi$  to a subset of $S$ of codimension $\geq 2$.

Without loss of generality we can assume also that $ M_ +$ and $ M_ -$ have
no common component.
 $L + M _+ \sim  D + M _-$.
For any curve $C$ in $S$ such that $\pi ^{-1} (C)$ is a smooth surface,
$L_ {|\pi ^{-1} (C) } = \sum \Sb s \endSb \sum \Sb k \endSb l _{k_s}$,  where $
l_{k_s}$ are the exceptional curves on the surface $ \pi ^{-1} (C)$
(with multiplicity) and $s\in C \cap     \bsa$.

For a fixed point $s$, $\pi^ {-1} (s) =  \sum \Sb k \endSb l _k \ + \ e$,
where  $e$ is the non-exceptional part of the fiber and the intersection matrix
of the $\{ l _k \}$'s
is negative definite (Zariski's Lemma [Z]).
Hence :
$$(L + M _+) \cdot \sum \Sb k \endSb l _k = L  \cdot \sum \Sb k \endSb l _k = (
\sum \Sb k \endSb l _k) ^2 < 0 \ .$$
Therefore there exists a  curve $ l _ {k_0} \in L$ which is in the base locus
of $L + M _+$; also
since the set of such $s$ is dense in $\pi (L)$, there exists a component of
$L$, say $L_1$,  in the fixed
component of $L + M _+$. Write
$$| L + M _+| = L _1 + \ |(L -L _1) \ + M _+| \ ;$$ then
the intersection matrix of the curves $  (L _ 1) _ {| \pi ^{-1} (s)}$ is
negative definite.
Proceeding by induction on the number of components, we conclude that $L$ is in
the fixed component
of $| L + M _+|$.

Hence $| L + M _+| = L + \ |M _+|$.
But  the  linear system $ |M _+|$ consists only of the fixed component $M _+$
($M _+$  maps to points via $\pi$),
therefore each component of $M _-$ has to be a component of $L$ or of $ M _+$,
 contradicting the assumption.

\medpagebreak

(1.13.2) The above argument shows also that if $E \neq \emptyset$, then $K_X$
is not $\pi$-nef. Using the relative version of Mori's theory we know that
 $X$ is birationally equivalent
 to a projective threefold
$\bar \pi : \bx  \rightarrow  S$ such that $\bx$ has only $\Bbb Q$-factorial
terminal singularities, and $K_{\bx }$ is $\bar \pi$-nef [Mo; 0.3.10].
Without loss of generality we can assume that $\mu : X \to \bx$ is a morphism.
 Then $K_ X\equiv \mu ^* (K_{\bx}) + \sum e_i E_i, \ e_i>0 \ \forall i$.
 Comparing with 1.13.1 we have
$$ \mu ^* (K_{\bx} - \bar \pi ^*(\kla) ) = E - \sum e_i E_i.$$
We need to prove that $E$ is exceptional for $\mu$.
This follows from [G2; proof of Theorem 1.1 (Case 2), pages 294-295].
\qed \enddemo

\bigpagebreak

\proclaim{Corollary 1.14}  $E-G$ (as in 1.7) is an effective $\Bbb Q$-divisor,
 if and only if  there exists a birationally equivalent fibration over
 $ \bar \pi : \bx \to S $ such that $\bx$ has only $\Bbb Q$-factorial
terminal singularities and   $K_{\bx} \equiv \bar \pi ^*  (K_ S +  \bla ) .$
\endproclaim

\smallpagebreak

 If $\pi : X \to S $ is an elliptic fibration between smooth varieties, with no
multiple fiber over divisors and $\bla$ with
simple normal crossings, then there exists
a birationally equivalent fibration $\bx \to S$ such that
 $K_{\bx} \equiv \bar \pi ^*  (K_ S +  \bla ) .$ In fact $E-G$ is effective.

\vskip 0.15in

\noindent (1.15) Following Miranda we refer to the the double points of the
ramification
locus with simple normal crossing as the ``\it collision points". \rm

The divisors of singular fibers on $X$ which map to irreducible components of
the ramification divisor intersecting at a point $p$, are said to \it``collide
at $p$". \rm
 In general the fiber over the intersection point of the ramification
 divisor is not one of Kodaira type.(Take a birational model of example 1.2.1,
with a simple normal crossing ramification divisor.)

\medpagebreak

The varieties in consideration are assumed to be complex projective, normal and
$\Bbb Q$-factorial (see [Wi]).

\smallpagebreak

\vskip  0.2in

\noindent { \bf \S 2. Blowing down: minimal model program and flat models.}

\bigpagebreak

The goal of this paragraph is  proving Theorem 2.3 and Corollary 2.7.
 We state the theorem in  general form, and then we will apply it to the case
studied
by Miranda in section 4.

\proclaim{Lemma 2.1, [N1, 0.4], [YPG; Hr, III, ex 10.9]}  Let  $\bar \pi :
\bx_0 \to \bs $ be an elliptic fibration such that $\bx_0$ has terminal
singularities.

 If $K_{\bx} \equiv \bar \pi ^*(D)$, for some $\Bbb Q$-divisor $D$ on $\bs $,
then $D \equiv K_{\bs } +  \bla _{\bx_0 / \bs } $ and
$(\bs , \bla _{\bx_0 / \bs })$ has log terminal singularities (as in 1.9.2).

If $\bs $ is smooth and $\bar \pi$ is equidimensional,
then $\bar \pi$ is also flat.
\endproclaim

\proclaim{Theorem 2.2}
Let $\bar \pi : \bx \rightarrow S$ be an elliptic fibration, where $\bx$ is a
threefold
with terminal singularities such that $K_{\bx} \equiv \bar \pi ^* ( \kla)$.

Let $\psi : S \rightarrow \so$ be a $\kla $-extremal contraction of an
irreducible curve $\g$.

(2.2.1) Then
there exists a collection of extremal rays  $\{ R_j \}$ on $\bx$, such that
${\bar \pi } _* (R_j ) = \g  $, for all $ j$.

(2.2.2) Let $\mu : \bx \dasharrow \bx_0$ be the contraction map induced by all
 the extremal rays mapping to $\psi _*(\g)$, where
 $\bx _0$ is a threefold with terminal singularities.
If the support of a divisor
$E$ dominates $\g$, $ E$ is exceptional for $\mu$.

(2.2.3) The induced fibration $\bar \pi _0 : \bx_0 \rightarrow \so$ is
$1$-dimensional over $\psi (\g)$ and $K_{\bx_0} \equiv \bar \pi _0 ^*(K_{\so} +
\bla _{\bx_0 /  \so}).$
\endproclaim

\demo{Proof}

(2.2.1) Let $\varepsilon = \psi \cdot \bar \pi$; any curve $\g _\a$ in
 $ \bar \pi ^{-1} ( \g  )$ is
mapped by $\varepsilon$ to a point; furthermore if $\g _ \a$
 is not a fiber of $\bar \pi $,
$$ K _{\bx} \cdot \g  _\a = \bar \pi  ^ *( K_ S +
\bla) \cdot \g _\a <  0.$$
Thus $K_{\bx}$ is not $\varepsilon$-nef and there exist a threefold $\bx_0$ and
a morphism $\pi _0$
as in the diagram
 $$
\CD
{\bx _0 }     @<\mu<<    \bx \\
@V{\bar \pi_0}VV           @V{\bar \pi }VV   \\
\so        @<\psi<<       S  \\
\endCD
$$
such that $K_{\bx_0}$ is $\bar \pi _0$-nef [Mo].

 $ \g _ \a$ belongs to the negative part of the cone $\overline {NE} (\bx)$ and
we
 can write $ \g  _\a = \sum \Sb i  \endSb q_i + \sum \Sb k \endSb R_k ^{\bx}$,
where $R_i ^{\bx}$ are extremal rays and $q_i$ other non extremal generators of
the cone.
Now,
$$ \g  = {\bar \pi}  _* \g  _\a = \sum \Sb j \endSb {\bar \pi} _* ( q _j) +
\sum \Sb k \endSb {\bar \pi }_* R_k ^{\bx},$$
since $\g$ is extremal the above equality implies that $\g$ and ${\bar \pi} _*
R_k ^{\bx}$ belong to the same extremal ray on $S$, for all $i$.
In particular $ \g  \sim c {\bar \pi} _* R_i ^{X_0}$, for some constant $c$ and
thus $ \g  =  \pi _* R_i ^{\bx}$.

\smallpagebreak

(2.2.2) Let $E$ be a divisor which dominates $\g $; without loss of generality
we can assume $E$ irreducible.
 We will assume
that $E$ is not contracted by $\mu$ and derive a contradiction.

 By hypothesis there exists an infinite collection of curves $\g _\a$ on $E$
such that $ K_{X_0} \cdot \g _\a \ < 0$. Note that $\mu$ is a composition of
divisorial contraction and flips [Mo] and denote by $ \mu _i: \bx _{ {(i-1)}}
\rightarrow \bx _{i } $  the first divisorial contraction, $E_i$ the
exceptional divisor and $\bar \pi _i$ the induced fibration.  There exists  an
infinite subcollection $\{ \g _\beta \}$
such that
 $ K_{X_{ {(i-1)}}} \cdot \g _\beta \ < 0$ [G2, Lemma   0.6] and it is thus
possible  to choose the $ \g _\beta$ so that no one of
them is contained in any of the $E _i$ (we have assumed that $E$ is not
contracted by $\mu$).
  Therefore:
$$ 0 >{ K_{\bx_{(i-1)}} \cdot \g _\beta } - { E _i \cdot \g _\beta} =
K _{\bx _i} { \mu _i } _* \cdot ( \g _\beta) , \ \forall\g _\beta .$$

\noindent Iterating this process we get:
$$0 > K_{ \bx _0} \cdot  \mu  _* (\g  _\gamma ), \text{ for some } \gamma.$$

As above, write  $ \mu  _*  \g  _\gamma = \sum \Sb i  \endSb g_i + \sum \Sb k
\endSb R_k ^{ \bx _0}$,
where $R_i ^{\bx _0}$ are extremal rays.
Then:
 $$ 0 =  \sum \Sb j \endSb \bar {(\pi _0 )}_* (g_j) + \sum \Sb k \endSb \bar {(
\pi _0 )}_* R_k ^{ \bx _0}.$$
Hence $ {(\bar \pi _0 )}_*  R_k=0$, and $R_k$ belongs to the fiber over $\psi
(\g)$  contradicting the definition of $\mu$.
Note that the extremal rays mapping to $\psi (\g)$ are all the extremal
rays mapping to points, because $K_{\bx} \equiv \bar \pi ^* ( \kla)$.

\smallpagebreak

(2.2.3)
By  2.2.2, $\mu$ contracts all divisors dominating $\g$.
 The idea of the proof is
showing that all the divisors in $ { \bar \pi }^ {-1} \g$ intersecting this
first set of exceptional divisors are contracted as well, and so on.

\smallpagebreak

Set:
 $$\aligned \Cal E _{\bx}  &= \{ \bar E \subset \bx , \text{  irreducible
divisors such that } \bar  E \in
  {\bar \pi } ^{-1} (\g)\} \\
\Cal S  _{ \bar E }  &= \{ \bar E ^k \in \Cal E _{\bx} \text{ such that
there exits a chain } (\bar E ^k, \cdots, \bar E  ^1 = \bar E) \text{ of
divisors}  \\   &\phantom{mn} \in \Cal E _{\bx}\text{ such that } \bar E ^j
\cap \bar E^ {j-1}
\neq \emptyset \text{ and } \bar E ^k \text{ dominates } \g \}.
\endaligned
$$
For any divisor $\bar E \in \Cal E _{\bx} $,  $ \Cal S _ {\bar E } \neq
\emptyset $,
 because the fibers  are connected. In fact, the connected component of
 $  \Cal E _{\bx} $containing $\bar E$ would consist only of divisors mapping
to points
 via $\bar \pi$.
 For any divisor $ \bar E \in \Cal E _{\bx} $ define the $alt$ of $\bar E$ as:
$$alt (\bar  E ) = min \{ M  \text{ such that } ( \bar E ^M, \dots, \bar E ^0 =
\bar E ) \text{ is a chain as above } \}. $$
Note that such $M$ is bounded.
\smallpagebreak

By (2.2.2) $\mu$ contracts all divisors dominating $\g$, i.e all divisors of
alt zero.
   Assume by induction that $\mu$ contracts all
 the divisors of alt $( M - 1) $ but not the ones of alt $M$, and derive a
contradiction.

 Let $\bar E$  be a divisor such that $alt(E)=M$ and  $ ( \bar E ^{M}, \dots
,\bar E ^0 =\bar E ) $ the corresponding minimal chain.  Then
$$K _{ \bx }  = [\mu ^* ( K _{\bx _0 } ) + \sum b _i  \bar E ^i ] \ ,$$
where $b _i > 0$; note that each $ \bar E^i \in \Cal S  _{ \bar E }$ is in the
summand and that by assumption $ \bar E^i \neq \bar E, \ \forall i$.

Note that there exists a family of curves $\{ C _{\gamma} \}$ in $\bar E$ such
that  $\bar E ^i \cdot  C _{\gamma} \geq  0 $, for
any $\bar E ^i $  in the summation and $\bar E ^j \cdot C _{\gamma} >  0 $
for some $\bar E ^j$ in the chain.
Note also that $K _{\bx} \cdot C _{\gamma}  = 0$, since $alt ( \bar E ) > 0$.
 By [G2, Lemma 1.6] (stated below) this implies that $ 0 > \mu ^*  ( K _{\bx ^0
} ) \cdot C _{\beta} = K _{\bx ^0 } \cdot \mu _* ( C _{\beta} )   $
 for any curve $ C _{\beta} $ in some subfamily
$\{ C _{\beta} \}$  of $\{ C _{\gamma} \}$, contradicting
  $ ( K _{ \bx _0 }  ) \text {   } \bar \pi _0$- nef.

 It follows also that
$K_{\bx_0} \equiv {\bar \pi _0}^*(K_{\so} + \bla _{X_0 / \so})$. \qed
\enddemo

\proclaim{ Lemma [G2, 1.6]} Let $\mu : X \rightarrow X ^+ $ be a ``flip"  map
between
threefolds.  If $\{ C _{\g} \}$ is an infinite collection of curves such that
 $ K  _X \cdot C _{\g}  < 0 $, for all $\g$, then there exists an infinite
subcollection
 $\{ C _{\beta} \}$  of $\{ C _{\g} \}$ such that for every $ C _{\beta} $,
 $ K _{X ^+} \cdot \mu _*  ( C _{\beta} )  < 0 .$
\endproclaim

\medpagebreak

If $\psi$ is a ($\kla$)-extremal contraction of $\g$, then in particular $\kla
\neq \psi ^*(K_{\so} + \bla _0)$ and $K_X \neq \epsilon ^*( K_{\so} + \bla
_0)$.

The following proposition is in some sense the converse of the previous one.

\proclaim{Proposition 2.3}
 Let $X$ and ${\bx}$ be threefolds with terminal singularities and
$$
\alignat3
& {X} && &&\  \ \ \ {{\bx}}\\
& \pi \searrow  && &&   \swarrow{\epsilon}  \\
&  &&\ \ \so && \\
\endalignat
$$
two  birationally equivalent elliptic fibrations.
Assume that  $K_{\bx}$  and
$K_X$ are numerically equivalent to the pullback of a
 $\Bbb Q$-divisor on $\so$.

Then $\epsilon$ is equidimensional if and only if $\pi$ is equidimensional
(flat if $\so$ is smooth).
\endproclaim

\demo{Proof} Proposition 2.3 is a particular case of the following Proposition
2.4.

 In fact $K_{\bx}\equiv {\epsilon}^*(K_{\so} +
\bla _{\so})$ and  $K_X \equiv \pi ^* (K_{\so} +
\bla _{\so})$, by 2.1 [N1, 0.6] and (1.9).\qed\enddemo

\proclaim{Proposition 2.4}
 Let ${\bx}$  be a variety with terminal singularities and
$$
\alignat3
& {X} && &&\  \ \ \ {{\bx}}\\
& \pi \searrow  && &&   \swarrow{\epsilon}  \\
&  &&\ \ \so && \\
\endalignat
$$
two birationally equivalent fibrations.
Assume that $K_{\bx} \equiv {\epsilon} ^ *(D)$,
$K_X \equiv \pi ^*(D)$,
for some $\Bbb Q$-divisors $D$ on $\so$  and that the fiber of $\epsilon$  over
a point $Q$
on $\so$ does not contain any divisor.
 Then the fiber of $\pi$ over a point $Q$
 does not contain any divisor.
\endproclaim

\demo{Proof}
\noindent
 Let $\tilde X$ be a common resolution of ${\bx}$ and $ X $  and let
$  \{ F _i \}$ be the exceptional irreducible divisors of the  morphism  $f:
\tilde X \to {\bx}  $  and $ \{ G _i  \}$  of  $g: \tilde X \to X $.

Then
$$\aligned
  K _{\tilde X} \equiv & f ^* ( K _{ X '} ) + \sum d _i F _i \equiv  f ^* \cdot
{\epsilon} ^*(D ) + \sum d _i F _i \\
\equiv  &g ^* ( K _{ X} ) + \sum c _i G _i  \equiv  f ^* \cdot {\epsilon} ^* (D
) + \sum c _i G _i
\endaligned
$$
If $T $ is an irreducible divisor $X$ mapping to the point $Q$,  then $T=F _i$,
for some index $i$.
Since  ${\bx}$ has terminal singularities, then $F_i$  has a positive
coefficient in the
 formula of the canonical bundle and has to be exceptional also for $g$.
Contradiction.
  \qed
\enddemo

\remark{Remark} The above proposition 2.4 has no assumption on the dimensions
and
the nature of the fibration.\endremark

\medpagebreak

\proclaim{Theorem 2.5}
 Let  $\bar \pi : \bx \rightarrow S$  an elliptic fibration of threefolds with
terminal singularities such that $K_{\bx}\equiv \bar \pi ^* (K_S + \bla _{\bx /
S})$.

 Let $ \psi :S  \to \so$ be a birational morphism with
irreducible exceptional divisor $\g$ and write:
$$ \kla _{\bx / S} = \psi ^ * (K_{\so} + \bla _{\bx / \so}) + \delta \g .$$

\smallpagebreak

(2.5.1)  If $\delta \geq 0 $
 there exists a birationally equivalent elliptic fibration $ \bar \pi _0 :
\bx_0 \rightarrow \so$, such that $\bx_0$ is a threefold with terminal
singularities and $$K_{\bx_0} \equiv \bar \pi ^*(K_{\so}
+ \bla _{\bx_0 / \so}).$$

Moreover,  $\delta > 0$ if and only if the fiber  over $\psi (\g)$ of any such
$\bar \pi _0$ is $1$-dimensional.

\smallpagebreak

(2.5.2) If $\delta  <0$, no birational equivalent fibration
 $\bar \pi _0 : \bx _0 \to \bs _0$ ($\bx _0$ with terminal singularities) can
satisfy $K_{\bx_0} \equiv \bar \pi _0 ^*(K_{\so}
+ \bla _{\bx_0 / \bs _0})$. Furthermore $\psi \cdot \bar \pi$ cannot have an
equidimensional model.
\endproclaim

\bigpagebreak

\remark{Remarks}

(2.6.1) In the above theorem $\delta > 0$ if and only if $\psi$ is a
($\kla$)-extremal contraction (1.13).

(2.6.2) We do not make assumptions on the equidimensionality of $\bar \pi$ over
$\g$.

(2.6.2) 2.5.1 always holds if $\so$ is smooth, $\bsa _0$ has normal crossing
and  the fibration has a  section (1.13), (4.3). 2.5.2 can actually occur, see
4.3.2 and
Example 1.11.
\endremark

\bigpagebreak

\demo{Proof} Set  $\epsilon = \psi \cdot \bar \pi : \bx \to \bs _0$.

\smallpagebreak

(2.5.1): If $\delta =0$, then
$$K_{\bx} \equiv \bar \pi ^*(\kla _{X / S}) \equiv \bar \pi ^* (\psi ^*(
K_{\so} + \bla _{X_0 / \so})) \equiv \epsilon  ^*( K_{\so} + \bla _{X_0 /
\so})$$
and the fiber of $\psi \cdot \bar \pi$ over $\psi(\g)$ contains a divisor.
 The statement follows from 2.3.
 Proposition 2.3 proves the other part of the statement.

\smallpagebreak

(2.5.2): Set $\beta = -b <0$. Then
 $K_{\bx} \equiv \epsilon ^* (K_{\so} + \bla _{X_0 / \so}) - b \bar \pi
^*(\g)$.
$K_{\bx}$ is $\bar \pi$-nef by hypothesis and $\epsilon$-nef.   All
birationally
equivalent fibration over $S_0$ have
a $2$-dimensional component over $\psi _*(\g)$, [Hm, 3.4].

 Assume that $\bar \pi _0 : \bx _0 \to \so$ is a birationally equivalent
fibration such that $\bx _0$ has terminal singularities and $K_{\bx_0} \equiv
\bar \pi _0 ^*(K_{\so}
+ \bla _{\bx_0 / \so})$. Let $X$ be a common resolution, as in the diagram:
$$
\alignat3
& &&\ \ X &&\\
 & ^{f}  \swarrow && &&\searrow ^{g}\\
& {\bx} &&  &&\  \ \ \ {\bx _0}\\
& _{\epsilon } \searrow && && \swarrow _{\bar \pi _0}\\
& &&\ \ \so &&
\endalignat
$$

Then $f ^*(K_{\bx}) = g ^*(K_{\bx _0})$ [Ko, 4.4]. That is, by commutativity of
the
diagram
$$ f ^* \cdot \epsilon ^*( K_{\so} + \bla _{X_0 / \so}) - b g ^* \bar \pi _0
^*(\g) \equiv
 f ^* \cdot \epsilon ^*( K_{\so} + \bla _{X_0 / \so}),$$
contradicting $b \neq 0$.
 \qed \enddemo

 Nakayama shows the existence
of an equidimensional model but with a \it different \rm base surface [N2].

\bigpagebreak

\proclaim{Corollary 2.7} Let $X$ be any elliptic threefold. Then then there
exists an equidimensional birational equivalent elliptic fibration $\bar \pi :
\bar  X_n \to   \bs _n$ such that  $\bx_n$ has terminal singularities,
$K_{\bx_n} \equiv
\bar \pi ^*(K_{  \bs _n} + \bla _{ \bs _n})$,  and
$(  \bs _n, \bla _{\bs _n})$ has log terminal singularities. Then, either

 (2.7.1) the Kodaira dimension of $X$ is non negative and $\bx_n$ is minimal.
Or

 (2.7.2) there exists  a morphism $\phi : \bs _n \to C$ such that $dim( S_n)<
dim(C)$
and the general fiber of $\phi \cdot \bar \pi$ is Fano.
If $\bs  _n$ is smooth, then $\bar \pi$ is flat.
\endproclaim

\demo{Proof} Nakayama [N2] proves that there always exists an equidimensional
birational equivalent fibration $\bx \to \bs $ satisfying the first part of the
statement.
We need to show that (2.7.1) and (2.7.2) hold for some particular model.

If $K_{\bx}$  is nef (that is, $K_{S} + \bla _{\bx / S}$ is nef), then $\bx$ is
minimal and we are done.
Otherwise we apply Mori's algorithm to $(S, (K_{\bs } + \bla _{\bs })$.
 Each birational contraction $\bs   \to \bs _1$ is thus log extremal and by
2.3.2 there exists an  elliptic
fibration $\bx_1 \to \bs _1$ satisfing the same part of the theorem. We can
repeat the process until     we obain a minimal elliptic threefold $\bx_n \to
S_n$  or  there exists  a morphism $\phi : \bs _n \to C$ such that $dim(\bs
_n)<  dim(C)$. In the second case $X$ is uniruled [Mo].
  The remaining statements follow from Lemma 2.1. \qed \enddemo

Corollary 2.7 is a stronger version of Theorem 1.1 [G2], but the proof
presented there is
different.  Nakayama's argument [N2] implies 2.7.1, but not 2.7.2.

\vskip 0.15in

\proclaim{Proposition 2.8}
 Let  $\bar \pi : \bx \rightarrow S$ be an elliptic fibration where $\bx$
is a threefolds with terminal singularities and $K_{\bx}\equiv \bar \pi ^* (K_S
+ \bla _{\bx / S}) +D$, with $D$ non effective.

 Let $ \psi :S  \to \so$ be a birational morphism with
irreducible exceptional divisor $\g$ and write:
$$ \kla _{\bx / S} = \psi ^ * (K_{\so} + \bla _{\bx / \so}) + \delta \g .$$

If $\delta  \leq 0$, no birational equivalent fibration
 $\bar \pi _0 : \bx _0 \to \bs _0$ ($\bx _0$ with terminal singularities) can
satisfy $K_{\bx_0} \equiv \bar \pi _0 ^*(K_{\so}
+ \bla _{\bx_0 / \bs _0})$. In particular $\bar \pi$ cannot have an
equidimensional
birationally equivalent fibration over $S_0$.

If $\delta >0$ such a  $\bar \pi _0 : \bx _0 \to \bs _0$ exists only if
$\delta \bar \pi ^*(\g) -G$ is effective.
\endproclaim

\demo{Proof} It follows from  1.14. Note that if an elliptic
fibration is equidimensional then $E-G$, as in 1.13, is effective.
\qed \enddemo

\vskip 0.2in

{\bf \S  3. Blowing up: Miranda's work. }

\smallpagebreak

 In this section  we review properties of Weierstrass models and Miranda's
construction:
given an elliptic fibration with section, it is possible to choose a
Weierstrass
model $W \to T$, where $W$ has rational singularities and $T$ is smooth.  The
fibration
$W \to T$ is thus flat, because it is equidimensional [Hr, III, ex 10.9].

 The starting point of Miranda's construction is blowing up $T$ to get a
 birationally equivalent fibration $\tilde W \to \so$ with
ramification locus with simple normal crossings, and then resolve the
singularities of the threefold $\tilde W$. Note that the pullback of a
Weiestrass model is a Weiestrass
model.

Over a general point of each component of  $ \bsa _{\tilde W/ \so}$
the singularities are uniform, and  so can be resolved \it  as surface
singularities, \rm
maintaining the flatness of the fibration.
Over the collision points the task is not as easy, but Miranda shows, checking
case by case, that after blowing up the surface a sufficient number of times it
is possible
to resolve the singularities of the new threefold preserving the flatness.

Let $\tilde \pi : \tilde X \to \so$ be the smooth  threefold thus obtained (the
 {\it smooth Miranda model).} Then
$$K_{\tilde X} \equiv  \tilde \pi ^* (K_\so + \bla _0), \tag 3.1$$
 in fact $G = E =\emptyset$ in formula 1.7, because the $\tilde \pi$ is flat
and the
singularities over the general points of $\bsa _{\tilde W/ \so} $ are resolved
as
 surface singularities.

We say that $\bx _0$ is a \it terminal Miranda model \rm if  $\bx _0 \to \so $
is a flat birationally equivalent elliptic fibration, satisfying
3.1 with  terminal singularities along the collision.  (Smooth points are
terminal
singularities, see for example [W].)

\smallpagebreak

Miranda writes explicitly a local analytic equation for each collision.  Using
the explicit equations, he shows  that the singularities of the threefold at
the collisions denoted by $\bullet$ in the table below can be resolved
preserving the flatness of the fibration.
 The next step in his algorithm is blowing up the double points of the
ramification divisor and describe the type of the singular fiber on (the
generic point of) the exceptional locus.

 The possible collisions  and the types of the singular fiber on the
exceptional locus are listed in the following table.
 Thus by further blowing up the surface all the collisions are reduced to the
one denoted by $\bullet$ in the table below.

\medpagebreak

$$
\vmatrix
  &  &  II & IV  & I_0^* & IV ^*   & II ^* & III  & I_0^*  & III^* &  I _c & I
^*_d   \\
&\\
  &  \\
II ^*     &  & I_0 & II & IV &  I_0 ^* &  IV^* &  & & & &\\
IV ^*    &   & \bullet II ^* & I_0 & II & IV &   &  &  & & &\\
I_0^*    &  & \bullet IV ^* & \bullet  II ^* & I _0 &  &  &  & & & & \\
IV   &  & \circ I _0 ^* & \circ  IV ^* &  &   &   &  &  & & &\\
  II  &  &  \circ  IV &  &  &  &   &  & & & & \\
III^*   &  &  &  &  &  &   & I_0& III & I_0 ^*  & \\
I_0 ^*   &  &   &  &  &  &   & \bullet III ^*& I_0 &  && \\
III      &  &   &  &  &  &   & \bullet I_0  ^* & &  & &\\
I ^*_b  &  &   &  &  &  &   &  & & &\bullet  I ^*_{b+c} &  I_{b+d} \\
I_a    &  &   &  &  &  &   &  & & & \bullet I_{a+c}  & \bullet I^*_{a+d}
 \endvmatrix
$$

It is easy to verify that also the collisions denoted by $\circ$, have a
Miranda model. Note that the only terminal (and not smooth) singularities at
the collisions occur at $II$-$II$ and $IV$-$IV$ (see Remark after Corollary
4.6);
flatness follows from equidimensionality also in this case, (see also 2.1).

\vskip  0.2in

{\bf \S 4.  Variations on an algorithm of Miranda.}

\bigpagebreak

We now apply Theorem 2.5 to the case of elliptic fibration $X_0 \to \so$
with section over a smooth surface and ramification divisor with simple normal
crossing.
If $\pi _0$ has a section, then  $n_i =0$ in (1.7) and
$\Cal O_{\so} (\bla _{\so}) =\Cal O_{\so} (\bda _{\so})$ is a line bundle .

\smallpagebreak

\proclaim{Hypothesis 4.0} All throughout \S 4, $\pi _0 : X_0 \rightarrow \so $
is an elliptic fibration between smooth varieties, $dim(X) =3$ and $ \bsa _0 $
is a divisor of simple normal crossings.
 $Q$ is a double point of $ \bsa _0 = \bsa _{X / \bs _0}$, and
 $\psi : S\rightarrow \so $ be the blow up of $\so$ at $Q$, with exceptional
divisor $\g$.

Let $X$ be the resolution of the pullback of $X_0$ by $\psi$, as in the
diagram:
$$
\CD
X_0    @<<<    X\\
@V\pi _0VV         @V \pi  VV      \\
\so     @<\psi<<          S
\endCD
$$

As in \S 1, we set $\bla = \bla _{X/S}$ and $\bla  _0= \bla _{X _0/\so}$.

It is clear that $\bla _0 = \psi _*(\bla)$ (1.9.4); in this section we compare
$\psi ^*(\bla _0)$ with $\bla$.  Write
$$\psi ^* (\bla _0) - \bla = \a \g, \text{ for some } \a \in \Bbb Q .$$
\endproclaim

\medpagebreak

\proclaim{Notation 4.1} Following 1.8, we write $\bla _0 = \sum a_k D_k
+ \frac{1}{12} J_{\infty}(\so) + \sum \frac{n_i -1}{n_i} Y_i $

and $$\bla = \sum a_k \hat D_k + a( \g ) \g + \frac{1}{12} J_\infty (S) +
\sum \frac{n_i -1}{n_i} \hat Y_i + \frac{n(\g) -1}{n(\g)} \g,$$ where  $\hat {
D_i}$ (  $\hat { Y_i}$) denotes
the strict transforms of the divisor $D_i$ ($Y_i$), $n(\g)$ the multiplicity of
the general fiber along $\g$ and  $a(\g)$ is the coefficient
associated to $\g$ (as in Table 1, \S 1). We can write
$$\psi ^* ( \bla _0 )  =  \sum a_k  \hat D _k +  \frac{1}{12} J_\infty (S)  +
\beta \g, $$
for some $ 0 < \beta \in \Bbb Q$; in fact $J_\infty (S) = \psi^* J_\infty
(\so)$, see \S 1.

Comparing the two formulas, we have
$$ \beta -a(\g) - \frac{n(\g) -1}{n(\g)} = \a .$$
\endproclaim

\medpagebreak

 \proclaim{Lemma 4.2} Suppose that $\pi _0$ has a section.
Let $a _ \ell, \ a_m$ be the coefficients determined by
the monodromy around
each of the two branches of the ramification divisor containing $Q$.
 Then $\beta = (a _ \ell + a_m)$, $\a = \lco \beta \rco$ and thus $a(\g) = \{
\beta \} =  \{ a_ \ell + a_m \} $, where $\{ x \}$  denotes the fractional part
of $x$ and $ \lco x \rco = x -\{ x \}$ .

 In particular
 $\a = 0$ or $1$ and  $\a \g $ is an effective divisor.
\endproclaim

\noindent {\bf Note:}  $\psi ^* (\bla _0) - \bla  = \a \g$ is in general only a
$\Bbb Q$-divisor.
 For example this is the case when the fibers over the generic points of the
components of the ramification divisor containing $Q$  are  multiple, or more
generally when the fibration does not admit a section (see 4.3).

\smallpagebreak

 The proof of the Lemma shows that the monodromy around the 2 branches of the
ramification divisor
containing $Q$ determine the type of monodromy around the exceptional divisor
$\g$.
 If there exists a section \it in a  neighborhood \rm of $Q$
the monodromy matrices and the type of singular fiber determine each other
uniquely (as we see in the proof of the Lemma below). This is not true in
general, see Proposition 4.3.

\smallpagebreak

\demo{Proof}

(4.2.0) Choose local coordinates such that $ Q = (0,0) $ and the two branches
of the ramification divisor  at $Q$ are defined by   $ {\bsa _0 }^ \ell : \{
x_0 = 0 \}$ and $ {\bsa _0 }^ m  =  \{y_0 = 0$ \}.
 Denote by $\gamma _{x_o }, \ \gamma _{y_o } $ , the monodromy matrices over
the generic point of $    {\bsa _0 }^ \ell $ and $    {\bsa _0 }^ m $.

The blow up of $\so$ at $Q$, in the coordinate chart $x_1 =x_o , \ y_1 =y_o  /
x_o $, induces the change in monodromy
$\gamma _{x_1}  = \gamma _{x_o } \cdot\gamma _{y_o }, \ \gamma _{y_1}  = \gamma
_{y_o }$.
Therefore the monodromy matrix around the generic point of $\g$,
the exceptional divisor of the blow up, is the composition of the two monodromy
matrices $\gamma _{x_o } $  and $\gamma _{y_o } $.

\smallpagebreak

\noindent Case (4.2.1): The monodromy matrix around the general point of $
{\bsa _0 }^ \ell $ (resp.  $    {\bsa _0 }^ m $) has finite order  ([Kd], 7.3
pg.  1281).

The only monodromy matrices in $SL(2, \Bbb Z)$ of finite order (2, 4, 3, 6) are
the ones listed in [Kd].
(It is enough to show it for $SL(2, \Bbb Z) \ / \ \pm I_d$, see [L, \S 4 Thm
1].)
Thus, the type of fiber, the monodromy matrix and its eigenvalue (in the closed
upper half plane) determine each other uniquely.  Furthermore, for each fixed
$\bold J$ value, any two monodromy matrices commute.
The eigenvalue of the matrix around the general point of $    {\bsa _0 }^ \ell
$  (resp.  $    {\bsa _0 }^ m $) is $e^ {2 \pi i a_ \ell}$ (resp.  $e^ {2 \pi i
a_m}$)  (see Table 1, \S 1).
Therefore the eingenvalue of the matrix around the general point of $\g$ is $e^
{2 \pi i (a_ \ell + a_m)}$, and the coefficient of
$\g$ in  2.2.3.  is $\{ a_ \ell + a_m \} $.

In the notation of 4.2, we have $a(\g) = \{ a_ \ell + a_m \} $,
 $\beta = (a _ \ell + a_m)$ and $\a = \lco a_ \ell + a_m \rco$.

Note that, since  $a _ \ell $ and $ a _m$ are rational numbers both  strictly
less then 1, then $\lco a _ \ell +  a _m\rco $ = 1 or 0 and $\psi ^* (\bla _0)
- \bla $ is an effective divisor.

\smallpagebreak

\noindent Case (4.2.2): $J$ has a pole at $Q$.

 Let $B_j$ be the irreducible divisor supported on ${\bsa _0} ^j$. Recall
that  $J_\infty (S) = \psi^* J_\infty (\so)$, see \S 1. By 4.2.0, the
monodromy
actions around  ${\bsa _0} ^\ell $ and ${\bsa _0} ^m $ determine the monodromy
action around $\g$ and thus $a(\g)$, by Table 1, in Section 1.

Case (i): The fibers over the general point of ${\bsa _0} ^j, \ j = \ell, m$ is
of type
 $I_{b_j}$. In this case $a(B_m) = a(B_\ell) =0$ and thus $\beta = a(\g) =0$.
Then $\psi ^* (\bda _0) = \bda .$

Case (ii) The fiber over the general point of  ${\bsa _0} ^\ell $ is of type
while the fiber over the
general point of ${\bsa _0} ^m $ is of type $I_{b_m} ^*$ . Then $a(B_m) = 1/2$,
$a(B_\ell) =0$ and  $\beta = a(\g) =1/2$.

Case (iii)  The fibers over the general point of ${\bsa _0} ^j, \ j = \ell, m$
is of type
$ I_{b_j} ^*, \ j = \ell, m$. In this case $a(B_m) = a(B_\ell) =1/2$ and thus
$\beta = 1$, while the combined monodromy actions give $a(\g) = 0$.
Consequently $\a = 1$.
 \qed
\enddemo

\bigpagebreak

More generally, the following hold,

\proclaim {Proposition 4.3} Notation as in 4.0.
Then $0 \leq \a < 2$, with $\a \in \Bbb Q$; actually, $0 \leq \a \leq 1$ except
in case 4.3.2,
below.
\endproclaim

\demo{Proof} Lemma 4.2. relates the behavior of the monodromy transformations.
This determines the type of the singular fiber over $\g$ if there is a section
around $Q$. In that case $\a = 0,1$.

\smallpagebreak

Case (4.3.1): The general fiber over the two branches of the ramification
divisor at $Q$ is not multiple, but the general fiber over $\g$ is multiple.

 Such examples can occur [DG], but  the collision has to be of type
$III$-$III^*$, $IV$-$IV^*$, $I_0^*$-$I_0^*$, because $a(\g)=0$, by Lemma 4.2
and Table 1, \S 1.  In this case  $\beta =1$ and $0 \leq \a = \frac{1}{n(\g)}
<1$.

Case (4.3.2): The general fibers over the two branches of the ramification
divisor are multiple, of multiplicity $n_1$ and $n_2$.

 Then $a_\ell = a_m = a(\g)=0$
 and $\beta = 2 - \frac{1}{n_1} - \frac{1}{n_2} >1$. Thus
 $$0 \leq \a = 2 -  \frac{1}{n_1} - \frac{1}{n_2} -1 + \frac{1}{n(\Gamma)}
<2.$$
The cases with $\a >1$ actually occur (see example 1.11); for a description
of the relations between $n_1$, $n_2$ and $n(\g)$, see [Gm].

Case (4.3.3): The general fiber over only one branch  (say $\Sigma _0 ^m$) is
multiple, of multiplicity $n_1$.

\noindent \phantom{Case} (i) Assume  $a_\ell >0$.
Then $n(\g) =0, \ a(\g) = a(\ell)$ and $$ 0\leq \a = a_\ell- a(\g) + 1
-\frac{1}{n_1} <1.$$

\noindent \phantom{Case} (ii) Assume $a_\ell =0$. Then $ a(\g) = a(\ell) =0$
 and $n(\g) = n_1$. Thus $\alpha =0$.
\qed \enddemo

\bigpagebreak

Note that  $\g $ is log extremal  if and only if
 $  (K _{S}  +  \bla ) \cdot \g <0$, that is
$$ \psi ^* ( K _{\so}  +\bla _0 ) \cdot \g +  ( 1 - \a ) \g  \cdot \g = ( 1 -
\a ) \g  \cdot \g  < 0  .$$
 Then  $1- \a >0$, that is $ \a <1$.

\smallpagebreak

\proclaim{ Corollary 4.4} With the notation of 4.0 and 4.2, assume that $\pi
_0$ has a section. The $\psi$ (the blow up of the collision point $Q$) is
($\kla$)-extremal if and only if $\a =0$. In particular

(4.4.1) If $ \bold J $ has a pole at $Q$, then $\psi$ is a log extremal
contraction unless the type of the fiber over the generic point of each
component of the ramification divisor containing $Q$ is of type $I_{b_j} ^*$,
for some $b_j$.

(4.4.1) If $ \bold J (  Q ) \neq \infty$, then $\psi$ is a log extremal
contraction if and only if $ \lco a_ \ell  + a_m \rco = \a = 0 .$
\endproclaim
\demo{Proof} $\psi$ is log extremal if and only if $0 \leq \a <1$, that is $\a
=0$, by Lemma 4.2 \qed \enddemo

\smallpagebreak

\proclaim{Corollary 4.5}
If $\pi _0$ has a section, $ \a  =0$ if and only if
 $X_0$ is birationally equivalent to an elliptic threefold
 $ \bar \pi _0 : \bx_0 \rightarrow \so$, with $1$-dimensional fiber over $\psi
(\g)$, such that  $\bx_0$ is a threefold with terminal singularities and
$K_{\bx_0} \equiv \bar \pi ^*(K_{\so}
+ \bla _{\bx_0 / \so})$.
\endproclaim

\demo{Proof} Since $X \to S$ does not have multiple fiber in codimension 1,
there exists a birationally equivalent fibration $\bar \pi : \bx  \to S $ such
that $K_{\bx} \equiv
\bar \pi ^* (\kla)$ (Theorem 1.13).
 By the above proposition 4.3, $\psi$ is a log extremal contraction if and only
if $\a=0$. Then  the statement follows from Theorem  2.5. \qed
\enddemo

\smallpagebreak

The following corollary shows which of the collisions have a flat (terminal)
resolution and which do not (the ``bad collisions").

\proclaim{Corollary 4.6}
Let $\pi _0: X _0 \to \so$ be a resolution of a Weierstrass model $W \to \so$
(as in \S 3).
Then:

(4.6.1) The collision above and on the diagonal in the following tables,
together
 with the collisions  of type  $I^*_{b_\ell}$-$I ^*_{b_m}$ cannot have a
Miranda model over $\so$.
(These are the ``bad" collisions.)

(4.6.2) The other collisions have a flat (over Q) terminal model $ \bar \pi _0
: \bx_0 \rightarrow \so$ such that  $K_{\bx_0} \equiv \bar \pi ^*(K_{\so}
+ \bla _{\bx_0 / \so})$.

(4.6.3) Blowing up a finite number of times it is possible to replace the
``bad" collisions,
with the other ones. The number of blowups needed depends only on
$ \lco a_ \ell + a _m\rco$.

\bigpagebreak

\noindent $$
\vmatrix
&\\
&   a _ \ell(type) &  && 1/4 (III) & 1/2 (I^*_0)  &3/4 ( III^* ) \\
&\\
a_m  (type) & \beta (type) & \\
& \\
&\\
3/4(III^*)  & & &  & 1(I_0) &  5/4 (III) & 3/2( I_0 ^* ) \\
&\\
1/2 (I_0^*) & & & &3/4 ( III^* )&  1(I_0) &   5/4 (III)\\
&\\
1/4(III )  & & &  &  1 / 2 (I_0^*) &3/4 ( III^* )&  1(I_0)
\endvmatrix
$$

\vskip 0.15in

\noindent$$
\vmatrix
&\\
& &  1/6 (II) && 1/3  (IV)  && 1/2 (I_0^*) && 2/3 (IV ^* )  &  5/6 (II ^*
)\\&\\
&\\
5/6(II ^*)  & &  1 (I_0) & & 7/6(II)  && 4/3(IV)  &&  3/2(I_0^*) &  5/3  (IV ^*
) \\
&\\
2/3(IV ^*) &  &  5 /6 (II^*) & & 1(I_0) & & 7/6 (II)& &4/3(IV) & 3/2(I_0^*) \\
&\\
1/2(I_0^*)  & & 2 /3 (IV^*)  && 5 /6 (II^*)  && 1(I_0)  && 7/6 (II) & 4/3(IV)
\\
&\\
1/3(IV)   & & 1 /2 (I_0^*)  && 2 /3 (IV^*)  && 5 /6 (II ^* )  && 1 (I_0)  & 7/6
\\
&\\
 1/6( II ) & & 1 /3 (IV) &  &1 /2(I_0)^*) &  &2 /3(IV^*)  & & 5 /6  (II ^* )&
1(I_0)
\endvmatrix
$$

\endproclaim

\medpagebreak

\demo{Proof} Recall that $\beta = (a _ \ell + a_m)$, $\a = \lco \beta \rco$ and
 $ X_0 \to \so$ has 1-dimensional fibers outside the collision points (see \S
3).

Let $Q$ be a collision point and $\psi$ the blow up at $Q$. The induced
fibration $ \pi :X \to S$ has also 1-dimensional fibers outside the collision
points
and $\g$.

\smallpagebreak

(4.6.1) In  case (4.6.1) $\psi$ is not a log extremal contraction (4.2-4.4). If
a Miranda model $ \tilde \pi : \tilde X \to \so$ exists, then  $\tilde \pi$
would be flat and
$K_{\tilde X} \equiv  \tilde \pi ^* (K_\so + \bla _0), $  contradicting 2.5.

(4.6.2) follows from Corollary 4.5.

(4.6.3)If $\alpha = \lco \beta \rco =1$,
then $a(\g) < max \{ a_\ell, a_m \} .$
\qed \enddemo

\vskip 0.1in

Mori's contraction algorithm outlined in 2.2 gives an explicit way to construct
a Miranda (terminal) model $\bar \pi : \bx _0 \to \so$ for all the ``good"
collisions.

 All the Miranda (terminal) models are isomorphic outside the point $Q$ and
related by flop transformations on the collision fiber over $Q$ (they are all
minimal models
relative to the respective fibrations) [Ko].

 In particular the flop transformations
preserve the type of analytic singularities [Ko].

\smallpagebreak

The Weierstrass model around the collisions $II$-$II$  is a terminal
 (but not smooth) Miranda model;  the resolution of the equisingular locus of
the collision of type $IV$-$IV$ is also a terminal (but not smooth) Miranda
model. All Miranda models (over the same base) of these collisions have to
 be terminal and not smooth [Ko].

\medpagebreak

All  the other ``good" collisions have smooth Miranda models.

\vskip 0.2in

{\bf \S   5. One more remark }

\medpagebreak

The equation of  example 1.2.1, $ y ^2 = x ^3 + sx + t$
 defines a smooth threefold $\pi _0 : X _0 \to \so$ in Weierstrass form.

 Kawamata's and Fujita's starting points (described in  \S 1)
use the hypothesis of the ramification divisor with simple normal
crossing. This is not the case for the above threefold, in fact there
exists a cusp at $s=t=0$.

 One can blow up $S_0$ three times to resolve the cusp,
  pull back the fibration, and at each step
  resolve the singularities of the threefold.
Let  $\psi _i : S_i \to S_{i-1}$, $i = 1,2,3$ be the blow ups.

 We then obtain  smooth threefolds $\pi _i : X_i \to S_i$,
$i = 1,2,3$.  It is easy to see that $\psi ^*(\bla _{i-1}) = \bla _i \ \forall
i$.

Each blow up $\psi _i : S_i \to S_{i-1}$ is then
in  a $K_{S_i} + \bla _i$-extremal contraction (1.10). Theorem 2.5 applies
and the Mori's algorithm used in 2.2 gives us back the original threefold
$X_0$. See [G1] for the explicit computations.

 \enddocument